\documentclass[12pt]{article}
\usepackage[utf8]{inputenc}

\usepackage{amssymb}

\renewcommand{\ge}{\geqslant}
\renewcommand{\le}{\leqslant}

{\clearpage\newpage\section*{Corrigendum -- submitted #1}}{}

\ifx\volno\undefined\def\volno{16}\fi
\ifx\volyear\undefined\def\volyear{2009}\fi
\ifx\papno\undefined\def\papno{R00}\fi

\newfont{\footsc}{cmcsc10 at 8truept}
\newfont{\footbf}{cmbx10 at 8truept}
\newfont{\footrm}{cmr10 at 10truept}

  \newlength{\BiblioSpacing}
  \renewenvironment{thebibliography}[1]{%
    \begin{oldthebibliography}{#1}%
      \setlength{\parskip}{\BiblioSpacing}
      \setlength{\itemsep}{\BiblioSpacing}
  }%
  {%
    \end{oldthebibliography}%
  }

\usepackage{amsthm,amsmath,amssymb}
\usepackage[colorlinks=true,citecolor=black,linkcolor=black,urlcolor=blue]{hyperref}

\theoremstyle{plain}
\newtheorem{theorem}{Theorem}
\newtheorem{lemma}[theorem]{Lemma}
\newtheorem{corollary}[theorem]{Corollary}

\newtheorem{conjecture}[theorem]{Conjecture}

\newtheorem{proposition}{Proposition}

\theoremstyle{definition}

\newtheorem{example}[theorem]{Example}

\theoremstyle{remark}

\title{\bf On the M\"obius Function of Permutations With One Descent}

\author{Jason P Smith\\
\small Department of Computer and Information Sciences\\[-0.8ex]
\small University of Strathclyde\\[-0.8ex] 
\small Glasgow, U.K.\\
\small\tt jason.p.smith@strath.ac.uk\\
}
\date{}

\begin{document}
	\maketitle

\begin{abstract}
The set of all permutations, ordered by pattern containment, is a poset. We give a formula for the M\"obius function of intervals $[1,\pi]$ in this poset, for any permutation $\pi$ with at most one descent.  We compute the M\"obius function as a function of the number and positions of pairs of consecutive letters in $\pi$ that are consecutive in value. As a result of this we show that the M\"obius function is unbounded on the poset of all permutations. We show that the M\"obius function is zero on any interval~$[1,\pi]$ where $\pi$ has a triple of consecutive letters whose values are consecutive and monotone. We also conjecture values of the M\"obius function on some other intervals of permutations with at most one descent. 
\end{abstract}

\section{Introduction}
Let $\sigma$ and $\pi$ be permutations of positive integers. We define an occurrence of $\sigma$ as a \emph{pattern} in $\pi$ to be a subsequence of $\pi$ with the same relative order of size as the letters in $\sigma$. For example, if $\sigma=213$ and $\pi=23514$ then there are two occurrences of $\sigma$ in $\pi$ as the subsequences $214$ and $314$. A permutation $\pi$ is said to \emph{avoid} a pattern $\sigma$ if there are no occurrences of $\sigma$ in $\pi$. The set of all permutations forms a poset $\mathcal{P}$, with a partial ordering defined as $\sigma\le\pi$ if $\sigma$ occurs as a pattern in $\pi$. An \emph{interval} $[\sigma,\pi]$ in $\mathcal{P}$ is a subposet consisting of all permutations $z\in \mathcal{P}$ with $\sigma\le z\le\pi$. The M\"obius function is defined recursively as follows: $\mu(\sigma,\lambda)=0$ if $\sigma\not\le\lambda$, $\mu(\sigma,\sigma)=1$ for all $\sigma$ and for $\sigma<\lambda$:
$$\mu(\sigma,\lambda)=-\sum_{\sigma\le z<\lambda}\mu(\sigma,z).$$
We frequently use the term \emph{M\"obius value} of a permutation $\lambda$ to refer to $\mu(1,\lambda)$ and we refer to permutations with a nonzero M\"obius value as \emph{nonzero} permutations. A \emph{descent} in a permutation $\pi=\pi_1\pi_2\ldots\pi_n$ is a decrease in the value of consecutive letters, that is, a descent at position $i$ is where~$\pi_i>\pi_{i+1}$.

Formulas for the M\"obius function in this poset in certain special cases have been proved. Almost all such results so far are on permutations constructed using \emph{direct sums}, where the direct sum of two permutations $\alpha$ and $\beta$, denoted $\alpha\oplus\beta$, is the concatenation of $\alpha$ with $\beta '$, where $\beta '$ is the permutation $\beta$ with each letter increased in value by the length of $\alpha$. For example $213\oplus 2413 = 2135746$. The first such result was by Sagan and Vatter in \cite{SagVat06}, where they give a formula for the M\"obius function on the poset of \emph{layered} permutations, that is, permutations that can be written as the direct sum of a number of decreasing permutations. More general results are presented in \cite{BJJS11} where a formula is given for the M\"obius function of all \emph{separable} permutations, that is, permutations avoiding both 3142 and 2413, along with many results for \emph{decomposable} permutations, that is, permutations that can be written non-trivially as direct sums. It is also shown that the absolute value of the M\"obius function has an upper bound in some of these cases. 
 
In this paper we present some of the first results for the M\"obius function on a substantial class of indecomposable permutations, the only other such result seems to be in~\cite{SteTen10}, which gives certain cases in which the M\"obius function is zero. As a result of this we show that $\mu(1,\pi)$ is unbounded, which does not seem to have been established before. Our main result is a formula for the M\"obius function on the interval $[1,\pi]$ for any permutation $\pi$ with at most one descent and that on such intervals the M\"obius function is alternating. Note that a permutation of length $n$ with one descent is indecomposable unless it starts with 1 or ends with $n$, so our result applies to a substantial class of indecomposable permutations.  

Define the subposet $\mathcal{P}_k\subseteq\mathcal{P}$ to be the poset containing permutations with exactly~$k$ descents. In this paper we mainly treat permutations from the subposets $\mathcal{P}_0$ and $\mathcal{P}_1$. We also use the notation $\mathcal{P}_k^n$ for the set of permutations of length $n$ with exactly $k$ descents and $\mathcal{P}^n$ for the set of all permutations of length~$n$.  As we often treat the M\"obius function of the intervals $[1,\pi]$ we  consider the M\"obius function as a function of a single variable in the form of $\mu(\pi):=\mu(1,\pi)$. An \emph{adjacency} in a permutation is two letters that are consecutive in position and have consecutive increasing values. For example, 24578136 has adjacencies~45 and~78, at positions 2 and 4. A permutation can also have a triple adjacency, or an adjacency of even greater length, as in the permutation 12456837 where there is the triple adjacency~456. The number of and positions of the adjacencies in a permutation will be key to our results. An important type of permutations from~$\mathcal{P}_1$ are the permutations without adjacencies, which are the permutations where odd and even letters are separated from each other. We denote the even length permutations without adjacencies as $M_n=246\ldots(2n)135\ldots(2n-1)$ and $W_n=135\ldots(2n-1)246\ldots(2n)$ for~$n>1$.

In Section \ref{sec:triple} we prove that $\mu(\pi)=0$ for any permutation $\pi$ containing a triple adjacency. In Section \ref{sec:noad} we prove a result relating to permutations with no adjacencies that is useful in the proof in Section \ref{sec:mob1}. In Section \ref{sec:mob1} we completely classify the M\"obius function on the intervals $[1,\pi]$ where $\pi$ has at most one descent. This proves the conjecture made in \cite{Ste13} that $\mu(\pi)=\binom{n+1}{2}$ when $\pi$ is of the form $246\ldots(2n)135\ldots(2n-1)$, which is the permutation without adjacencies $M_n$. This shows that $\mu(\pi)$ is unbounded in general, answering a question asked in \cite{BJJS11}, where it was shown that $|\mu(\pi)|\le1$ for all separable permutations $\pi$. In Section \ref{sec:furth} we present additional conjectures we have not been able to prove on the M\"obius function of permutations with at most one descent.

\section{The M\"obius function on permutations with a triple adjacency}\label{sec:triple}
In this section we present and prove a lemma stating that a permutation with a triple adjacency has a M\"obius value of zero. While interesting in its own right, it is useful in proving the result in Section \ref{sec:mob1}. But first we consider some notation and important points about adjacencies.

We defined an adjacency in the introduction as two letters that are consecutive in position and have consecutive increasing values. There is an analogous decreasing adjacency but we consider adjacencies to be increasing unless otherwise stated, because in~$\mathcal{P}_1$ decreasing adjacencies are rare and do not play a role in our considerations. We denote the value of an adjacency by the value of its initial letter, so in the permutation 24578136 the adjacencies 45 and 78 have values 4 and 7. Notice that a triple adjacency consists of two adjacencies of two letters, for example we can split 456 into 45 and 56. When counting the adjacencies in a permutation we count adjacencies of two letters, therefore 12456837 has three adjacencies 12, 45 and~56.

\begin{lemma}\label{lem:triple}
If a permutation $\pi$ contains a triple adjacency then~$\mu(\pi)=0$.
\begin{proof}
We can easily check that $\mu(123)=0$. Now assume that the claim holds for any permutation of length $m<n$ where $m\ge3$.
Given a permutation $\pi\in \mathcal{P}^n$ with a triple adjacency, removing any of the letters of the triple adjacency from $\pi$ results in the same permutation, call this $\sigma$. Any permutation obtained from $\pi$ by removing any of the letters not in the triple adjacency still has a triple adjacency hence by the inductive hypothesis has a zero M\"obius value. Hence all nonzero permutations in $[1,\pi)$ occur in $[1,\sigma]$, implying: $$\displaystyle\mu(\pi)=-\sum_{1\le z<\pi}\mu(1,z)=-\sum_{1\le z\le\sigma}\mu(1,z)=0.$$
\end{proof}
\end{lemma}
The result in Lemma \ref{lem:triple} also holds for the case of decreasing triple adjacencies, with an analogous proof. We can slightly generalise this result to give the following corollary, whose proof is analogous to the proof of Lemma \ref{lem:triple} after  suitably modifying the base case:
\begin{corollary}
If a permutation $\pi$ contains an adjacency (increasing or decreasing) of length $k\ge3$, then $\mu(12\ldots (k-2),\pi)=0$.
\end{corollary}

\section{The permutations with one descent and no adjacencies}\label{sec:noad}
We present a result on permutations with no adjacencies that is useful in proving the results in Section \ref{sec:mob1}. Before stating the lemma, we introduce some notation and definitions along with a few remarks about the posets $\mathcal{P}_0$ and $\mathcal{P}_1$.

We say that two permutations are \textit{related} if both or neither permutation begins with~1. For example the permutations $246135$ and $2357146$ are related as neither begins with 1 but $246135$ and $135246$ are not related.

Let the increasing permutation $12\ldots k$ be denoted $\mathbf{k}$. Notice that the poset $\mathcal{P}_0$ forms a chain, as for any $k\ge1$ the only length $k$ permutation without a descent is the increasing permutation $\textbf{k}$. As $\mathcal{P}_0$ is a chain it is easy to see that $\mu(\mathbf{k})=0$ for any $k>2$.

 It is also important to note that a permutation with $k$ descents cannot contain, as a pattern, a permutation with more than $k$ descents. Therefore, in any interval $[\textbf{1},\pi]$, with~$\pi\in\mathcal{P}_1$, any permutation $\lambda\in[\mathbf{1},\pi]$ must be in $\mathcal{P}_0\cup\mathcal{P}_1$. That is, the set $\mathcal{P}_0\cup\mathcal{P}_1$ is an order ideal in $\mathcal{P}$, also called a \emph{permutation class}. A permutation class can be uniquely determined by its \emph{basis}, that is, the set of minimal permutations it avoids. The basis for~$\mathcal{P}_0\cup\mathcal{P}_1$ can be shown to be  $\{321,2143,3142\}$. We remark that the poset of permutations with at most $k$ descents, for any fixed $k$, is a permutation class, but the basis for the general case $k>1$ is much more difficult to find. A formula is given in~\cite[Theorem~4.2]{BouFer13} which can be used to calculate the size of such a basis but this formula is rather complicated.

Recall that the permutations without adjacencies are the permutations where the odd and even letters are separated from each other. The even length permutations without adjacencies are $M_n=246\ldots(2n)135\ldots(2n-1)$ and $W_n=135\ldots(2n-1)246\ldots(2n)$ for~$n>1$.
\begin{lemma}\label{lem:noad}
Let $\pi\in P_1^n$ be a permutation with no adjacencies. Then $\pi$ contains, as patterns, precisely all permutations in $\mathcal{P}_1$ of length less than $n$ with at most two adjacencies except the following:
\begin{enumerate}
\item The permutations of length $n-1$ with two adjacencies.\label{noad1}
\item The permutations of length $n-1$ with one adjacency that are not related to~$\pi$.\label{noad2}
\item The permutations of length $n-2$ with two adjacencies that are not related to~$\pi$.\label{noad3}
\end{enumerate}
\begin{proof}
Let $R$ and $N$ be the subposets of $\mathcal{P}_1$ which contain the permutations, of length $m<n$, that are related and not related to $\pi$, respectively. Also denote the subposets of $R$ and $N$ with exactly $k$ adjacencies as $R^k$ and $N^k$, respectively. We need to prove that $\pi$ contains all permutations $\sigma\in R^0\cup N^0\cup R^1$, all permutations $\sigma\in R^2\cup N^1$ of length~$m<n-1$ and all permutations $\sigma\in N^2$ of length $m<n-2$.

First consider the permutations in $R$. Note that $R^0$ is a chain and any permutation in~$R^0$, of length $m<n$, can be obtained by removing the $n-m$ largest letters of $\pi$. To obtain a permutation $\sigma\in R^1$, of length $m<n$, where the adjacency has value $i$, it is necessary and sufficient to remove the letter $i+1$ from $\pi$ and then to adjust to the correct length permutation by removing the $n-m-1$ largest letters. So to create any permutation in $R^0\cup R^1$ there is only one letter that must be removed and thus all permutations in~$R^0\cup R^1$ of length $m\le n-1$ can be obtained from $\pi$. Now consider a permutation $\tau\in R^2$, of length $m<n$, with adjacencies of value $i$ and $k$. To create such a permutation, from $\pi$, we remove the letters of value $i+1$ and $k+1$ and then we adjust the length by removing the $n-m-2$ largest letters. So the permutations in $R^2$ require at least two letters to be removed and therefore all the permutations in $R^2$ of length $m\le n-2$ can be obtained, but none of length $n-1$.

Now consider the permutations in $N$. Removing the letter 1 from $\pi$ creates a unique length $n-1$ permutation $\lambda$ which is in $N^0$. We can then apply the same argument as in the previous paragraph to $\lambda$ instead of $\pi$. Hence we can obtain all permutations in~$N^0\cup N^1$ of length $m \le n-2$ and all permutations in $N^2$ of length $m\le n-3$. As $\lambda$ is the only permutation of length $n-1$ in $N^0$ we can get all permutations in $N^0$ of length $m \le n-1$.
 
Therefore we can obtain all permutations with at most two adjacencies in $R\cup N$ except for the following: The permutations in $R^2$ of length $n-1$, the permutations in $N^2$ of lengths $n-1$ and $n-2$ and the permutations in $N^1$ of length $n-1$.
\end{proof}
\end{lemma}

We provide an example of Lemma \ref{lem:noad}:
\begin{example}
Consider the permutation $\pi=135246$. By Lemma \ref{lem:noad} we know the only permutations with at most two adjacencies not contained in $\pi$ are:
\begin{enumerate}
\item The permutations of length 5 with two adjacencies that is: 12354,  41235, 12534, 34125, 12453, 31245,15234, 23415, 14523, 23145, 13452 and 21345.
\item The permutations of length 5 with one adjacency that are not related to~$\pi$ that is: 35124, 23514, 25134 and 24513.
\item The permutations of length 4 with two adjacencies that are not related to~$\pi$ that is: 4123, 3412 and 2341.
\end{enumerate}
\end{example}

\section{The M\"obius function for permutations with one descent}\label{sec:mob1}
In this section we present a formula for the M\"obius function on the interval $[\mathbf{1},\pi]$ where~$\pi$ is any permutation with at most one descent.
\begin{theorem}\label{thm:main}
Given a permutation $\pi\in \mathcal{P}_0\cup\mathcal{P}_1$, of length $n>2$, the value of $\mu(\pi)$ can be computed from the number and positions of adjacencies in $\pi$, as follows:
\begin{enumerate}
\item If $\pi$ begins with $12$ or ends in  $(n-1)n$ then $\mu(\pi) = 0$.\label{main1}
\item If $\pi$ has a triple adjacency then $\mu(\pi) = 0$.\label{main2}
\item If $\pi$ has more than two adjacencies then $\mu(\pi) = 0$.\label{main3}
\item If $\pi$ has exactly two adjacencies then:\label{main4}
\begin{enumerate}
\item If the first adjacency has greater value than the second then $\mu(\pi)=\pm1$,\label{main4a}
\item If the first adjacency has lower value than the second then $\mu(\pi)=0$.\label{main4b}
\end{enumerate}
\item If $\pi$ has exactly one adjacency, at position $i\in\{1,\ldots,n-1\}$, and the descent is at position $d$, then: (see item \ref{main7} for calculating the sign)\label{main5}
\begin{enumerate}
\item If $i<d$ and $\pi_1\not=1$ then $\mu(\pi)=\pm i$,\label{main5a}
\item If $i<d$ and $\pi_1=1$ then $\mu(\pi)=\pm (i-1)$,\label{main5b}
\item If $i>d$ and $\pi_{n}\not=n$ then $\mu(\pi)=\pm(n-i)$,\label{main5c}
\item If $i>d$ and $\pi_{n}=n$ then $\mu(\pi)=\pm(n-i-1)$.\label{main5d}
\end{enumerate}
\item If $\pi$ has no adjacencies then:\label{main6}
\begin{enumerate}
\item If $n$ is even and $\pi_1=1$, that is $\pi=W_{\frac{n}{2}}$, then $\mu(\pi)=-\dbinom{\frac{n}{2}}{2}$,\label{main6a}
\item If $n$ is even and $\pi_1=2$, that is $\pi=M_{\frac{n}{2}}$, then $\mu(\pi)=-\dbinom{\frac{n}{2}+1}{2}$,\label{main6b}
\item If $n$ is odd then $\mu(\pi)=\dbinom{\frac{n+1}{2}}{2}$.\label{main6c}
\end{enumerate}
\item If $\mu(\pi)\not=0$ then $\mu(\pi)$ is positive if and only if $n$ is odd.\label{main7}
\end{enumerate}
\end{theorem}
Before proving Theorem \ref{thm:main} we make some remarks:
\begin{itemize}
\item Each permutation with one descent falls into at least one of the above classes.
\item The above result agrees on permutations covered by more than one class. These cases are:
\begin{itemize}
\item A permutation with one adjacency and beginning with 12 or ending with $(n-1)n$ has zero M\"obius value by both part \ref{main1} and part \ref{main5}.
\item A permutation with 12 at the beginning and $(n-1)n$ at the end has zero M\"obius value by both part \ref{main1} and part \ref{main4b}.
\item A triple adjacency can be treated as two consecutive adjacencies and the result states the M\"obius value is zero by part \ref{main2} and part \ref{main4b}.
\item It is possible for a permutation to fall into all three of the first cases, such as 12354, and such a permutation has zero M\"obius value according to all three cases.
\end{itemize}
\end{itemize}

For part \ref{main1} of Theorem \ref{thm:main}, a permutation that begins with 12 or ends with $(n-1)n$ is decomposable so the proof follows directly from Corollary 3 in \cite{BJJS11} and part \ref{main2} follows from Lemma \ref{lem:triple}.

We prove the remaining parts of Theorem \ref{thm:main} using an inductive argument throughout the following subsections. For a base case we need to consider all permutations of length $3\le n \le 6$. We know certain permutations have zero M\"obius value by the already proven parts \ref{main1} and \ref{main2} of Theorem \ref{thm:main} So we can leave such permutations. We now list the remaining permutations of length $3\le n \le 6$ with one descent along with their calculated M\"obius value and which case of Theorem \ref{thm:main} they fall into:
$\mu(34125)=\mu(14523)=1  (\ref{main4a})$, $\mu(3412)=\mu(145236)=\mu(256134)=\mu(346125)=\mu(356124)=-1 (\ref{main4a})$, $\mu(235614)=\mu(236145)=\mu(361245)=0 (\ref{main4b})$, $\mu(231)=\mu(312)=\mu(13425)=\mu(14235)=\mu(23514)=\mu(25134)=1 (\ref{main5})$, $\mu(1423)=\mu(3124)=\mu(1342)=\mu(2314)=\mu(134625)=\mu(136245)=\mu(235146)=\mu(251346)=-1 (\ref{main5})$, $\mu(24513)=\mu(35124)=2 (\ref{main5})$, $\mu(245136)=\mu(351246)=\mu(146235)=\mu(135624)=-2 (\ref{main5})$, $\mu(132)=\mu(213)=1 (\ref{main6})$, $\mu(1324)=\mu(2413)=-1 (\ref{main6})$, $\mu(13524)=\mu(24135)=3 (\ref{main6})$, $\mu(135246)=-3 (\ref{main6a})$, $\mu(246135)=-6 (\ref{main6b})$.

It is straightforward to check that these results agree with Theorem \ref{thm:main}.

The reason it is necessary to check the base case up to length $n=6$ is so we can use Lemma \ref{lem:sumzero} to cancel out the M\"obius values of sets of permutations in the intervals we consider.

From now on we assume that any permutation in $\mathcal{P}_0\cup\mathcal{P}_1$, of length less than $n$, where $n>6$, satisfies the claims in Theorem \ref{thm:main} and we prove that Theorem~\ref{thm:main} then holds for permutations with at most one descent of length $n$ and thus for any length. When referencing the induction hypothesis we add the part of Theorem \ref{thm:main} being referenced in brackets, for example (\ref{thm:main}.\ref{main6a}) for Theorem \ref{thm:main} part \ref{main6a}.

By our inductive hypothesis (\ref{thm:main}.\ref{main3}) we can see that any nonzero permutation of length $m<n$ can have at most two adjacencies. If we combine this with Lemma \ref{lem:noad} we see that a permutation of length $m<n$ with no adjacencies contains all nonzero permutations of length at most $m-3$.

\subsection{The structure of the proof}
The remaining parts of the proof of Theorem \ref{thm:main} all follow a similar schema, which is outlined as follows:
\begin{enumerate}
\item Consider $\pi\in \mathcal{P}_1^n$.
\item Remove one letter from each adjacency in $\pi$ or the largest letter if $\pi$ has no adjacencies. This leaves a permutation $\lambda$ with no adjacencies.
\item By the definition of the M\"obius function, $\displaystyle\sum_{\sigma\in[\mathbf{1},\lambda]}\mu(\sigma)=0$.
\item Now we can compute $\mu(\pi)$ using $\displaystyle \mu(\pi)=-\sum_{\substack{\sigma\le\pi\\ \sigma\not\le\lambda}}\mu(\sigma)$.
\end{enumerate}

We develop this schema in detail for the proof of Proposition \ref{prop:main3} in the following subsection and then, as they are quite similar, the remaining parts of the proof are done in less detail.

We present two lemmas which we frequently reference throughout the proof:

\begin{lemma}\label{lem:altsign}
Let  $\sigma\in\mathcal{P}_1^m$, where $m<n$,  be a nonzero permutation satisfying either one of the following conditions:
\begin{enumerate}
\item Has exactly one adjacency, which is neither 12 nor $(m-1)m$.
\item Has exactly two adjacencies at least one of which is neither 12 nor $(m-1)m$.
\end{enumerate}
Then $\sigma$ contains a length $m-1$ permutation $\lambda$ with the same number of adjacencies as $\sigma$ such that $\mu(\lambda)+\mu(\sigma) = 0$.
\begin{proof}
If $\sigma$ has exactly one adjacency, at location $i$, then either this adjacency is before or after the descent. If the adjacency is before the descent, then by the induction hypothesis~(\ref{thm:main}.\ref{main5}) the M\"obius value of $\sigma$ is a function of $i$. We know $m$ must be to the right of $i$ so removing $m$ creates a length $m-1$ permutation $\lambda$ with exactly one adjacency at location~$i$, so by the induction hypothesis (\ref{thm:main}.\ref{main5} and \ref{thm:main}.\ref{main7}) $\mu(\sigma)=-\mu(\lambda)$. If the adjacency is after the descent then removing the letter 1 gives an analogous argument. This completes the first case.

If $\sigma$ has exactly two adjacencies we can remove either the letter 1 or $m$ which gives a length $m-1$ permutation $\lambda$ which has two adjacencies of the same relative sizes as the adjacencies in $\sigma$. By the induction hypothesis (\ref{thm:main}.\ref{main7}) the sign of the M\"obius function is alternating, therefore $\mu(\lambda)=-\mu(\sigma)$.
\end{proof}
\end{lemma}

\begin{example}
\begin{enumerate}
\item Consider the permutation $13425$ which is of the first form in Lemma~\ref{lem:altsign}. Removing the letter 5 gives the permutation $1342$. We compute the M\"obius values of these permutations as $\mu(13425)=1$ and $\mu(1342)=-1$.
\item Consider the permutation 24781356 which is of the second form in Lemma \ref{lem:altsign}. Removing the letter 1 gives the permutation 1367245. We compute the M\"obius values of these permutations as $\mu(24781356)=-1$ and $\mu(1367245)=~1$.
\end{enumerate}
\end{example}

We can use Lemma \ref{lem:altsign} to show that the M\"obius values of certain sets of permutations sum to zero.

\begin{lemma}\label{lem:sumzero}
Take a set $\Delta^1$ of $k$ nonzero permutations from $\mathcal{P}_1^m$, where $4<m<n$, all with~$t>0$ adjacencies and where none of the adjacencies is 12 or $(m-1)m$. Then we can construct the following sets:
\begin{itemize}
\item A set $\Delta^2$ of $2k$ permutations from $\mathcal{P}_1^{m-1}$ with exactly $t$ adjacencies.
\item A set $\Delta^3$ of $k$ permutations from $\mathcal{P}_1^{m-2}$ with exactly $t$ adjacencies.
\end{itemize}
Also the sum of the M\"obius values of all the permutations in $\Delta^1\cup\Delta^2\cup\Delta^3$ is zero.
\begin{proof}
From each permutation in $\Delta^1$ we have two options: We can remove the letter 1 or the letter $m$. Assume first that the removal of either of these letters from any of the permutations does not remove the descent from the permutation, then it is easy to see that this does not create or remove an adjacency. So to create $\Delta^2$ we get two permutations for each permutation in $\Delta^1$ by removing either~$1$ or $m$. To create $\Delta^3$ we remove, from each permutation in $\Delta^1$, both $1$ and $m$. It is easy to see that, as we are only removing the letters 1 and $m$, the permutations in the union of~$\Delta^2$ and $\Delta^3$ are distinct. This concludes the proof of the first part of the lemma.

To show that the M\"obius values sum to zero we can apply Lemma \ref{lem:altsign}. Recall that a permutation with $t>2$ adjacencies has M\"obius value zero by the induction hypothesis~(\ref{thm:main}.\ref{main3}). As $\Delta^1$ only contains nonzero permutations any permutation~$\lambda\in\Delta^1$ must be of one of the forms in Lemma~\ref{lem:altsign}. First suppose it is of the first form, that is, it has one adjacency, and suppose this adjacency is before the descent. We can pair $\lambda$ with the permutation obtained by removing the letter $m$ from $\lambda$ and their M\"obius values sum to zero. Then, given the permutation $\lambda^1\in\Delta^2$ obtained by removing the letter 1 from $\lambda$, we can pair this with the permutation $\lambda^{1,m}\in\Delta^3$ obtained by removing 1 and $m$ from $\lambda$. By Lemma~\ref{lem:altsign} we know the M\"obius values of these two permutations sum to zero. We can do this for each permutation in $\Delta^1$, which completes this case. An analogous argument applies to the case where $\lambda$ has an adjacency after the descent or has two adjacencies.

If the removal of the letter 1 or $m$ results in the removal of the descent from one of the permutations then we apply an analogous argument to entire proof above. In this argument we must account for the fact that for each permutation of this form there are two permutations that are increasing permutations of length greater than~$2$. As such, these permutations contain a triple adjacency and will have zero M\"obius value and the M\"obius value of the remaining permutations cancel as above.
\end{proof}
\end{lemma}

\subsection{Theorem \ref{thm:main} part \ref{main3}}
Recall that we are assuming Theorem \ref{thm:main} is true for any permutations of length $m<n$. We now consider part \ref{main3} of Theorem \ref{thm:main} for permutations of length $n$.
\begin{proposition}\label{prop:main3}
A permutation $\pi\in \mathcal{P}_1^n$ with more than two adjacencies has $\mu(\pi)=0$.
\begin{proof}
Suppose $\pi$ has $k>2$ adjacencies also suppose none of the adjacencies are $12$ or~$(n-1)n$ and there are no triple adjacencies. Then, by the inductive hypothesis (\ref{thm:main}.\ref{main3}), $\pi$ contains no nonzero permutations of length greater than $n-k+2$. There is a unique length $n-k$ permutation $\lambda$ contained in~$\pi$ with no adjacencies. Let us ignore $\lambda$ along with any other permutation in $[\mathbf{1},\lambda]$, since their contributions to the M\"obius value of $\pi$ sum to zero. Then we can use Lemma \ref{lem:noad} to consider the remaining permutations, that are possibly nonzero, occurring in $\pi$:
\begin{itemize}
\item Of length $n-k+2$ there remain $s=\binom{k}{2}$ permutations with two adjacencies, call these $\Gamma^0=\{\gamma_1^0,\ldots,\gamma_s^0\}$, where each of the $\gamma_i^0$'s is obtained by removing a letter from all but two of the adjacencies in $\pi$.
\item Of length $n-k+1$ there remain:
\begin{itemize}
\item $k$ nonzero permutations with one adjacency, call these $\Delta^0=\{\delta_1^0,\ldots,\delta_k^0\}$.
\item $2\binom{k}{2}$ permutations with two adjacencies obtained by removing the letter 1 or the largest letter from each of the $\gamma_i^0$'s, call these $\Gamma^1=\{\gamma_1^1,\ldots,\gamma_{2s}^1\}$. \end{itemize}
\item Of length $n-k$ there remain:
\begin{itemize}
\item All the permutations related to $\pi$ that have two adjacencies, where at least one of the adjacencies is an original adjacency in $\pi$, call this set of permutations~$\Omega^0$.
\item $\binom{k}{2}$ permutations not related to $\pi$ that have two adjacencies, both occurring in~$\pi$, call these $\Gamma^2=\{\gamma_1^2,\ldots,\gamma_s^2\}$.
\item $2k$ permutations with one adjacency obtained by removing the letter $1$ or the largest letter from the $\delta_i^0$'s, call these $\Delta^1=\{\delta_1^1,\ldots,\delta_{2k}^1\}$.
\end{itemize}
\item Of length $n-k-1$ there remain:
\begin{itemize}\item All permutations with two adjacencies, where at least one of the adjacencies is an original adjacency in $\pi$, call this set of permutations $\Omega^1$.
\item $k$ permutations with one adjacency obtained by removing the letter 1 and the largest letter from each $\delta_i^0$, call these $\Delta^2=\{\delta_1^2,\ldots,\delta_{k}^2\}$.
\end{itemize}
\item Of length $n-k-2$ all permutations not related to $\pi$ with two adjacencies, where at least one of the adjacencies is an original adjacency in $\pi$, call this set of permutations~$\Omega^2$.
\end{itemize}
Note that by Lemma \ref{lem:sumzero} the M\"obius values in $\Delta^0\cup\Delta^1\cup\Delta^2$ sum to zero and the same is true of $\Gamma^0\cup\Gamma^1\cup\Gamma^2$ and $\Omega^0\cup\Omega^1\cup\Omega^2$. We know these sets satisfy the length condtions in Lemma~\ref{lem:sumzero} because the maximum number of adjacencies is $\frac{n-2}{2}$, this is because the letters 1 and $n$ are not in adjacencies and there are no triple adjacencies, which implies $n-k\ge\frac{n+2}{2}\ge4.5>4$. This implies $\mu(\pi)=0$ and completes this case.

Now suppose one of the adjacencies in $\pi$ is $12$ or $(n-1)n$. If these adjacencies occur at the beginning or end, respectively, then this reduces to part \ref{main1} of Theorem \ref{thm:main}. It is also possible that one of these adjacencies occurs directly after or before the descent in which case the proof follows from the proof above with minor modifications. These modifications arise from the fact that removing the letter 1 from $\pi$ is equivalent to removing the adjacency 12 and likewise with the letter $n$ and the adjacency $(n-1)n$. In certain cases, this may result in $n-k\not>4$, and we must apply Lemma \ref{lem:altsign} to get the desired cancellation.
\end{proof}
\end{proposition}

\subsection{Theorem \ref{thm:main} part \ref{main4}}
\begin{proposition}\label{prop:main4}
Consider a permutation $\pi\in \mathcal{P}_1^n$ with exactly two adjacencies, at positions $k$ and $i$. If the first adjacency has greater value than the second then $\mu(\pi)=\pm1$, otherwise~$\mu(\pi)=0$.
\begin{proof}
If $\pi$ begins with 12 or ends with $(n-1)n$, then $\mu(\pi)=0$ by part \ref{main1} of Theorem \ref{thm:main}.  Now consider the case $\pi$ does not contain both the adjacencies 12 and $(n-1)n$.

Removing the letters $\pi_i$ and $\pi_k$ results in a permutation $\lambda$, of length $n-2$, with no adjacencies. As the M\"obius values of all the permutations in $[\mathbf{1},\lambda]$ sums to zero we can ignore any permutation in said interval. Now use  Lemma \ref{lem:noad} and consider the remaining permutations. By Lemma \ref{lem:sumzero} the M\"obius values of the remaining permutations with one adjacency sum to zero. Split the remaining permutations with two adjacencies into two sets $A$ and~$B$, where $A$ are those obtained from $\pi$ by removing the letters 1 or $n$ (or both) and $B$ are those obtained from $\pi$ by removing a letter from an adjacency and then removing another letter to create a new adjacency that does not occur in $\pi$. As the largest permutations in $B$ are of length $n-2>4$ we can apply Lemma~\ref{lem:sumzero} to see that the M\"obius values of the permutations in $B$ sum to zero. 

This just leaves us to consider~$A$. First assume $\pi$ doesn't have the adjacencies $12$ or~$(n-1)n$ directly after or before the descent. Then $A$ contains the following permutations with two adjacencies:
\begin{itemize}
\item A permutation $\delta$ of length $n-1$, obtained by removing the letter 1 from $\pi$.
\item A permutation $\tau$ of length $n-1$, obtained by removing the letter $n$ from $\pi$.
\item A permutation $\sigma$ of length $n-2$, obtained by removing letters $1$ and $n$ from $\pi$.
\end{itemize}
By the inductive hypothesis (\ref{thm:main}.\ref{main4} and \ref{thm:main}.\ref{main7}) it is clear that $\mu(\tau)+\mu(\sigma)=0$. This means that $\mu(\pi)=-\mu(\delta)$. The relative values of the adjacencies in $\pi$ are the same as in $\delta$ so, if the first adjacency has greater value than the second then $\mu(\pi)=-\mu(\delta)=\pm1$, otherwise $\mu(\pi)=-\mu(\delta)=0$. This completes the first case.

Now consider the case when $\pi$ contains the adjacency 12 but not $(n-1)n$, then $A$ only contains $\tau$ and $\mu(\pi)=-\mu(\tau)$. Similarly when $\pi$ contains the adjacency $(n-1)n$ but not 12, then $A$ only contains $\delta$ and $\mu(\pi)=-\mu(\delta)$. The result then follows by evaluating the value of $\delta$ or $\tau$. This completes this case.

Finally consider the case $\pi$ contains both adjacencies 12 and $(n-1)n$ and with $(n-1)n$ occurring before 12, in this case there are no permutations in the set denoted $A$ above and not all the permutations with one adjacency cancel. So we repeat the argument above considering the permutations with one adjacency.
\end{proof}
\end{proposition}
\textbf{Remark:} Note that the M\"obius value in the above proof is computed as a negation of a permutation of length one less. Hence the M\"obius value is alternating in the case of permutations with two adjacencies.

\subsection{Theorem \ref{thm:main} part \ref{main5}}
\begin{proposition}\label{prop:main5}
Consider a permutation $\pi\in \mathcal{P}_1^n$ which has exactly one adjacency at position $i$ and the descent at position $d$. Then:
\begin{enumerate}
\item If $i<d$ and $\pi_1\not=1$ then $\mu(\pi)=\pm i$,
\item If $i<d$ and $\pi_1=1$ then $\mu(\pi)=\pm (i-1)$,
\item If $i>d$ and $\pi_{n}\not=n$ then $\mu(\pi)=\pm(n-i)$,
\item If $i>d$ and $\pi_{n}=n$ then $\mu(\pi)=\pm(n-i-1)$.
\end{enumerate}
\begin{proof}
If $\pi$ begins with 12 or ends with $(n-1)n$, then $\mu(\pi)=0$ by part \ref{main1} of Theorem \ref{thm:main}. Next, we consider the case  where $\pi$ doesn't have the adjacencies $12$ or~$(n-1)n$ directly before or after the descent.

Removing $\pi_i$ from $\pi$ creates a permutation $\lambda$ with no adjacencies and we can ignore the interval $[\mathbf{1},\lambda]$ as the M\"obius values sum to zero by definition. We can apply Lemma \ref{lem:sumzero} to the remaining permutations with two adjacencies to see that their M\"obius values sum to zero. By  Lemma \ref{lem:noad} this leaves us to consider three permutations with one adjacency:
\begin{itemize}
\item Of length $n-1$ there remain two permutations with one adjacency, obtained by removing the letters 1 or $n$, call these $\sigma_1$ and $\sigma_2$ respectively.
\item Of length $n-2$ there remains one permutation with one adjacency not related to~$\pi$. This is obtained by removing the letters 1 and $n$ from $\pi$, call this $\delta$.
\end{itemize}
We consider the four cases in the statement of the proposition and obtain the M\"obius value from the induction hypothesis (\ref{thm:main}.\ref{main5}):
\begin{itemize}
 \item If $i<d$ then $\mu(\sigma_1)+\mu(\delta)=0$. Hence $\mu(\pi)=-\mu(\sigma_2)$ which gives:
 \begin{enumerate}
 \item If $\pi_1\not=1$ then $\mu(\pi)=-\mu(\sigma_2)=\pm i$.
 \item If $\pi_1=1$ then $\mu(\pi)=-\mu(\sigma_2)=\pm (i-1)$.
 \end{enumerate}
 \item If $i>d$ then $\mu(\sigma_2)+\mu(\delta)=0$. Hence $\mu(\pi)=-\mu(\sigma_1)$ which gives:
 \begin{enumerate}
 \setcounter{enumi}{2}
 \item If $\pi_n\not=n$ then $\mu(\pi)=-\mu(\sigma_1)=\pm (n-i)$.
 \item If $\pi_n=n$ then $\mu(\pi)=-\mu(\sigma_1)=\pm(n-i-1)$.
 \end{enumerate}
 \end{itemize}
 This completes this case of the proof.
 
If $\pi$ contains the adjacency 12 or $(n-1)n$ then removing the letter $n$ or 1, respectively, gives a permutation with one adjacency $\alpha$. The M\"obius values of all the other permutations sum to zero by Lemma \ref{lem:sumzero}, so $\mu(\pi)=-\mu(\alpha)$. Evaluating the four different cases of the proposition and using the inductive hypothesis (\ref{thm:main}.\ref{main5}) to get $\mu(\alpha)$ completes the proof.
\end{proof}
\end{proposition}
\textbf{Remark:} Note that in the above proof for each case the M\"obius value of $\pi$ is a negation of a permutation of length one less. Therefore the sign of the M\"obius value is alternating for all permutations with exactly one adjacency.

\subsection{Theorem \ref{thm:main} part \ref{main6}}
\begin{proposition}\label{prop:main6}
Let $\pi$ be a permutation in $\mathcal{P}_1^n$ with no adjacencies. Then:
\begin{enumerate}
\item If $n$ is even and $\pi_1=1$, that is $\pi=W_{\frac{n}{2}}$, then $\mu(\pi)=-\dbinom{\frac{n}{2}}{2}$,
\item If $n$ is even and $\pi_1=2$, that is $\pi=M_{\frac{n}{2}}$, then $\mu(\pi)=-\dbinom{\frac{n}{2}+1}{2}$,
\item If $n$ is odd then $\mu(\pi)=\dbinom{\frac{n+1}{2}}{2}$.
\end{enumerate}
\begin{proof}
First note that $\pi$ contains a permutation $\lambda$, with no adjacencies, of length~$n-~1$, obtained by removing the largest letter from $\pi$. As the M\"obius values of all the permutations in $[\mathbf{1},\lambda]$ sum to zero we can ignore any permutation in said interval. By  Lemma \ref{lem:noad} this leaves us to consider the following permutations which occur in $\pi$:
\begin{itemize}
\item Of length $n-1$ there remain:
\begin{itemize}
\item One permutation with no adjacencies obtained by removing the letter 1 from~$\pi$.
\item The permutations with one adjacency each obtained by removing a letter from~$\pi$, excluding the letters 1 and $n$.
\end{itemize}
\item Of length $n-2$ there remain:
\begin{itemize}
\item The permutations not related to $\pi$ with one adjacency. These are obtained by removing the letter 1 from each of the permutations with one adjacency of length~$n-1$ listed above.
\item All permutations of length $n-2$ related to $\pi$ with two adjacencies.
\end{itemize}
\item Of length $n-3$ there remain the permutations not related to $\pi$ with two adjacencies.
\end{itemize}

First consider the case $n$ is even and $\pi_1=1$, that is when $\pi=13\ldots(n-1)24\ldots n=W_{\frac{n}{2}}$. We will consider the permutations in $[1,W_{\frac{n}{2}}]$ based on number of adjacencies, and when needed, by the number removed to create an adjacency. We start with the nonzero permutations with two adjacencies. Note that all the length $n-3$ permutations with two adjacencies are obtained from the length $n-2$ permutations with two adjacencies by removing the letter~1. We can then apply Lemma \ref{lem:altsign} to see that the M\"obius values of the permutations with two adjacencies sum to zero. We can repeat this argument with the permutations with one adjacency obtained from $\pi$ by removing any of the letters $3,4,\ldots,(n-1)$ to see that these cancel with the permutations of length $n-2$ with one adjacency. This leaves the permutations obtained from $\pi$ by removing the letters 2 and 1, respectively. The first is of the form $124\ldots(n-2)3\ldots(n-1)$ and begins with 12 so has zero M\"obius value by part \ref{main1} of Theorem \ref{thm:main}. The second is $24\ldots(n-2)13\ldots(n-1)$ which by the induction hypothesis (\ref{thm:main}.\ref{main6c}) has M\"obius value $\binom{\frac{n-1+1}{2}}{2}$, which implies $\mu(W_{\frac{n}{2}})=-\binom{\frac{n}{2}}{2}$.

In the case $n$ is odd and $\pi_1=1$, the argument is analogous. We find the permutation~$24\ldots n13\ldots(n-1)$ has M\"obius value $-\binom{\frac{n-1}{2}+1}{2}$, which implies $\mu(\pi)=\binom{n+1}{2}$.

Next consider the case where $n$ is even and $\pi_1=2$, that is $\pi=24\ldots n13...(n-1)=M_{\frac{n}{2}}$. First we consider permutations of length $n-1$ with one adjacency formed by removing one of the letters $3,4,\ldots,(n-1)$. We can apply Lemma \ref{lem:altsign} to see that the M\"obius value of all but one of these cancel with all but one of the length $n-2$ permutations with one adjacency. The only remaining length $n-2$ permutation is $124\ldots(n-2)35\ldots(n-3)$ which has zero M\"obius value by part \ref{main1} of Theorem \ref{thm:main}. The only remaining length $n-1$ permutation  is $24\ldots(n-2)(n-1)13\ldots(n-3)$ which by the induction hypothesis~(\ref{thm:main}.\ref{main5a}) has M\"obius value $\frac{n}{2}-1$.

Now consider the remaining permutations with two adjacencies. The permutation with the triple adjacency 123 contributes zero to the M\"obius value by part \ref{main1} of Theorem~\ref{thm:main}. Removing the letter 2 and any letter $i>3$ from $\pi$ results in a permutation with adjacency~12 immediately after the descent. If $i$ is even then the larger adjacency also appears after the descent so such a permutation contributes zero to the M\"obius value by the induction hypothesis (\ref{thm:main}.\ref{main4b}). If $i$ is odd then the adjacency appears before the descent. Since each such permutation has M\"obius value $-1$ by the induction hypothesis~(\ref{thm:main}.\ref{main4a}), and~$i$ is an odd number between 5 and $n$, the sum of the M\"obius values of these permutations is $-\frac{n}{2}+2$. We can apply Lemma \ref{lem:altsign} to cancel all of the other permutations with two adjacencies in a similar way to the case $\pi_1=1$ above. 

We must also consider the M\"obius values of the permutations found by removing 2 or~1 from $\pi$. The permutation $35\ldots(n-1)124\ldots(n-2)$ has M\"obius value $\frac{n}{2}-1$ by the induction hypothesis (\ref{thm:main}.\ref{main5c}). The permutation $13\ldots(n-1)24\ldots(n-2)$ has M\"obius value~$\binom{\frac{n}{2}}{2}$ by the induction hypothesis (\ref{thm:main}.\ref{main6c}). The M\"obius value of $\pi$ is given by the negation of the sum of the M\"obius values of the permutations it contains, so we sum the above values and negate which gives:
$$\mu(\pi)=-\left(\dbinom{\frac{n}{2}}{2}+2(\frac{n}{2}-1)-\frac{n}{2}+2\right)=-\dbinom{\frac{n}{2}+1}{2}.$$

Finally we consider the case where $n$ is odd and $\pi_1=2$, that is, $\pi=24\ldots(n-1)13\ldots n$. The argument proceeds in an analogous manner to the previous case, except the sum of the M\"obius values of the permutations with two adjacencies is $\frac{n-1}{2}-2$ and the M\"obius values of the three permutations with one adjacency are $-\frac{n}{2}+1$, $-\frac{n}{2}+1$ and~$\binom{\frac{n-1}{2}}{2}$, resulting in~$\mu(\pi)=\binom{\frac{n+1}{2}}{2}$.
\end{proof}
\end{proposition}
\textbf{Remark:} \begin{itemize}
\item The nice form of the result in Proposition \ref{prop:main6} raises the question of a direct combinatorial proof. We expect to present such a proof in the forthcoming paper \cite{Smi14} which analyses topological properties of some intervals in the poset $\mathcal{P}$.
\item Notice that in the above cases the M\"obius value is positive if and only if $n$ is odd. Therefore the M\"obius value is alternating.
\end{itemize}

\subsection{Finishing the proof of Theorem \ref{thm:main}}
Notice that the remarks after Propositions \ref{prop:main4}, \ref{prop:main5} and \ref{prop:main6} show that the M\"obius value is alternating for all nonzero permutations, which implies the M\"obius value is positive if and only if $n$ is odd. This proves part \ref{main7} of Theorem \ref{thm:main}. We have shown that if the classification of Theorem \ref{thm:main} holds for all permutations, of length less than $n$, with at most one descent, then it also holds for $n$. By induction, that completes the proof of Theorem~\ref{thm:main}. Parts \ref{main5} and \ref{main6} of Theorem \ref{thm:main} give us the following important corollary:
\begin{corollary}
On the poset $\mathcal{P}$ the function $\mu(\pi)$ is unbounded.
\end{corollary}

\section{Conjectures on the M\"obius function for intervals of permutations with at most one descent}\label{sec:furth}
So far we have mainly concentrated on intervals of the form $[\textbf{1},\pi]$. We now consider permutations where we allow the permutation $\textbf{1}$ to change. We see that this change increases the complexity of computing the M\"obius function quite drastically especially in the second conjecture we present, but also leads to some interesting results relating to the M\"obius function being dependant on whether a permutation is separable.
\subsection{The M\"obius function on the intervals  \texorpdfstring{$[\sigma,M_n]$}{Lg} and  \texorpdfstring{$[\sigma,W_n]$}{Lg}}\label{sec:mobbot}
In this subsection we examine intervals $[\sigma,\pi]$ where $\pi$ is one of the two permutations of \emph{even} length with no adjacencies and $\sigma\in\mathcal{P}_1$. Recall that these permutations with no adjacencies are denoted $M_n=24\ldots(2n)13\ldots(2n-1)$ and $W_n=13\ldots(2n-1)24\ldots(2n)$. This leads us to the following conjecture which has been checked by computer to hold for any pair $(m,n)$ where $m<12$ and $n<7$:
\begin{conjecture}\label{con:bottom}
  Given a permutation $\sigma\in \mathcal{P}_1^m$, let $i$ be the number of adjacencies in $\sigma$.  If $\sigma\le\pi$ where $\pi\in\{M_n,W_n\}$ we have the following:
\begin{itemize}
\item If $\sigma$ is separable, then:
\begin{itemize}\renewcommand{\labelitemii}{$\diamond$}
\item $\mu(\sigma,M_n)=\pm\dbinom{n+1}{m}$,
\item $\mu(\sigma,W_n)=\pm\dbinom{n+m-i-2}{m}$.
\end{itemize}
\item If $\sigma$ is not separable:
$$\mu(\sigma,\pi)=\pm\dbinom{n+\lfloor\frac{m-i-a}{2}\rfloor}{m}$$
where
$a=\begin{cases}
0,&\mbox{ if } \sigma \mbox{ and } \pi \mbox{ are related}\\
1,&\mbox{ otherwise }
\end{cases}.$
\end{itemize}
Also the M\"obius value is positive if and only if $m$ is even.
\end{conjecture}
Recall that when considering an adjacency of length $k$ we regard it as $k-1$ individual adjacencies.

Notice that Conjecture \ref{con:bottom} only deals with intervals $[\sigma,\pi]$ where $\pi$ is of even length. In Theorem \ref{thm:main} we can see that changing $\pi$ between odd and even length has little effect on the M\"obius function.  In Conjecture \ref{con:bottom}, on the other hand, there is a substantial difference  between the odd and even case.

\subsection{The M\"obius function on the interval \texorpdfstring{$[M_m,\pi]$}{Lg}}\label{sec:mobtop}
We can reverse the idea in subsection \ref{sec:mobbot} and consider intervals $[\sigma,\pi]$ where $\sigma$ is a permutation without adjacencies and $\pi\in\mathcal{P}_0^n\cup\mathcal{P}_1^n$.  In this subsection we conjecture a formula for the M\"obius function on such
intervals.  This formula is somewhat complicated, but turns out to be computationally efficient, compared to the brute force method of computing from the recursive formula for the M\"obius function. Before stating the result we define a few statistics on $\pi$:
\begin{itemize}
\item Let $a$ be the number of adjacencies in $\pi$.

\item Set $\hat{n} =\begin{cases}n-1,&\mbox{ if }\pi_{n}=n\\
n,&\mbox{ otherwise}\end{cases}$.

\item Let the set $A=\{i_1,\ldots,i_a\}$ be the ordered sequence of the values of the adjacencies in $\pi$. Also add to $A$ two phantom adjacencies $i_0$ and $i_{a+1}$ which occur before and after the descent, respectively, with values:\newline  $i_0 =\begin{cases}-1,&\mbox{ if }\pi_1\not=1\\
0,&\mbox{ otherwise}\end{cases}$ and $i_{a+1} = \hat{n}+1$.

\item A function: $$\widehat{C}_\beta^\alpha(k,s)=\begin{cases}
\dbinom{\alpha-2k}{\beta},&\mbox{ if } 0\le k<\frac{s}{2}\\
\dbinom{\alpha-2(s-k)+1}{\beta},&\mbox{ if } \frac{s}{2}\le k<s\end{cases}.$$

\item A sequence $\hat{J}=\{\hat{j}_0,\ldots,\hat{j}_{a}\}$ where:
 $$\hat{j_k}=\left\lfloor\frac{i_{k+1}-i_k-2}{2}\right\rfloor.$$

\item Split $\hat{J}$ into two sequences $j^a$ and $j^b$ in the following way:
$$\begin{cases}\hat{j_k}\in j^a,&\mbox{ if } i_{k} \mbox{ and } i_{k+1} \mbox{ occur on the same side of the descent}\\
\hat{j_k}\in j^b,&\mbox{ otherwise}\end{cases}.$$

\item Set $\displaystyle s=\sum_{t=0}^{a}\hat{j_t}$.

\item Trim $j^a$ and $j^b$ in the following way:
\begin{enumerate}
\item If $j^a$ is empty remove the largest element from $j^b$ and set $\epsilon=0$,
\item If $j^a$ is not empty let $max_{j^a}$ be the largest element in $j^a$ and remove it from~$j^a$, then set $\displaystyle \epsilon = max_{j^a}-\sum_t j^a_t$,
\item Then set $\alpha=|j^a|$ and $\beta=|j^b|$,
\item Remove all zero elements from both sequences and if this results in $j^b$ being empty set $\epsilon =0$,
\item Finally sort $j^a$ into ascending order and $j^b$ into decreasing order.
\end{enumerate}

\item  Define the function $\displaystyle s_\theta(\kappa,\tau)=\sum_{t=\kappa}^{\tau}j^\theta_t$.

\item Set $\lambda=\left\lceil\frac{\hat{n}}{2}\right\rceil+m-\left\lceil\frac{5a}{2}\right\rceil+\beta-t$ and $\sigma = 2m-2a+\beta$\newline
where $t=\begin{cases}1,&\mbox{if } \pi_1=1 \mbox{ and } n \mbox{ is even and } \pi_n=n\\
1,&\mbox{if } \pi_1=1 \mbox{ and } n \mbox{ is odd and } \pi_n\not=n\\
0,&\mbox{otherwise}\end{cases}$.
\end{itemize}
For an example of these statistics see Example \ref{ex:stat} below.
 
We can now state the conjecture which has been checked by computer tests to hold for all pairs $(m,n)$ where $m<6$ and $n <12$:
\begin{conjecture}\label{con:top}
Consider the interval $[M_m,\pi]$ where $\pi\in\mathcal{P}_0^n\cup\mathcal{P}_1^n$ and $\lambda$, $\sigma$, $j^a$, $j^b$, $\epsilon$, $s$, $s_\theta$ and $\widehat{C}$ are all as defined above, then:

If $\pi$ begins with 12, ends with $(n-1)n$ or contains a triple adjacency $\mu(M_m,\pi)=0$, otherwise:
\begin{align*}
|\mu(M_m,\pi)|&=\dbinom{\lambda}{\sigma}-
\sum_{\tau=0}^{|j^b|-1}
\sum_{\gamma=0}^{\tau}
\sum_{\omega=\tau-\gamma}^{j^b_\gamma+s_b(\tau+1,|j^b|-1)-1}
\widehat{C}^{\lambda-\tau-2}_{\sigma-\tau-1}(\omega,s)\\
&+\sum_{\tau=0}^{|j^a|-1}\left[
\sum_{\gamma=1}^{j^a_{\tau}+s_a(0,\tau-1)}\widehat{C}^{\lambda-|j^b|-\tau}_{\sigma-|j^b|-\tau}(\gamma,s+1)
+\sum_{\omega=1}^{\epsilon}\widehat{C}^{\lambda-|j^b|-\tau}_{\sigma-|j^b|-\tau}(\omega+1,s+1)
\right].
\end{align*}
Also the sign of $\mu$ is positive if and only if $n$ is even.
\end{conjecture}

Whilst Conjecture \ref{con:top} is rather complicated it is significantly more computationally efficient than computing the M\"obius function from its recursive definition. To see this consider the following example, for an interval of rank 20, whose computation from the recursive definition would take enormous time even on a fast computer:
\begin{example}\label{ex:stat}
  Consider the  interval 
$$I=[24681357,2\,4\,6\,7\,9\,12\,14\,16\,18\,21\,23\,24\,26\,28\,1\,3\,5\,8\,10\,11\,13\,15\,17\,19\,20\,22\,25\,27].$$
 We compute $\mu(I)$ using Conjecture \ref{con:top}, first extracting the following statistics from $I$:
\begin{itemize}
\item $a=4$, $m=4$ and $\hat{n}=28$,
\item $A=\{-1,6,10,19,23,29\}$ and $\hat{J}=\{2,1,3,1,2\}$,
\item Before trimming: $j^a=\{2,3\}$ and $j^b=\{1,1,2\}$,
\item After trimming: $j^a=\{2\}$, $j^b=\{2,1,1\}$, $\epsilon=1$, $\alpha=1$ and $\beta=3$,
\item $s=9$, $\lambda=11$ and $\sigma=3$.
\end{itemize}
Putting this into the formula of Conjecture \ref{con:top} we get:
\begin{align*}
\mu(I)&=\binom{11}{3}-\sum_{\tau=0}^{2}\sum_{\gamma=0}^{\tau}\sum_{\omega=\tau-\gamma}^{j^b_\gamma+s_b(\tau+1,|j^b|-1)-1}\widehat{C}^{9-\tau}_{2-\tau}(\omega,9)\\
&+\sum_{\tau=0}^{0}\left[\sum_{\gamma=1}^{j^a_{\tau}+s_a(0,\tau-1)}\widehat{C}^{8-\tau}_{0-\tau}(\gamma,10)+\sum_{\omega=1}^{1}\widehat{C}^{8-\tau}_{0-\tau}(\omega+1,10)\right]\\
&=\binom{11}{3}-\sum_{\omega=0}^{4}\widehat{C}^9_2(\omega,9)-\left[\sum_{\omega=1}^{2}\widehat{C}^8_1(\omega,9)+\sum_{\omega=0}^{1}\widehat{C}^8_1(\omega,9)\right]\\
&-\sum_{\omega=0}^{0}\widehat{C}^7_0(\omega,9)+\sum_{\gamma=1}^{2}\widehat{C}^{8}_{0}(\gamma,10)+\sum_{\omega=1}^{1}\widehat{C}^{8}_{0}(\omega+1,10)\\
&=\binom{11}{3}-\binom{9}{2}-\binom{7}{2}-\binom{5}{2}-\binom{3}{2}\\ &-\binom{6}{1}-\binom{4}{1}-\binom{8}{1}-\binom{6}{1}-\binom{7}{0}+\binom{6}{0}+\binom{4}{0}+\binom{4}{0}\\
&=165-36-21-10-3-6-4-8-6-1+1+1+1=73.
\end{align*}
\end{example}
Whilst we cannot verify this is the correct value of the M\"obius function on this interval, the example serves as a good indicator of the efficiency of the conjecture if it can be proved correct.

\section*{Acknowledgements}
I would like to thank Eva Hauksd\'ottir for suggesting some of the results in the statement of Theorem \ref{thm:main} and my advisor Einar Steingr\'imsson for many helpful comments. I would also like to express my gratitude to an anonymous referee for extremely detailed, useful and in depth comments, for suggesting the rewriting of the proof of Proposition 4 that is presented here and pointing out an error, now fixed, in Lemma 8.

\end{document}